\newif
\ifconfver
\confvertrue        

\ifconfver
    \documentclass[10pt,journal,final,twoside]{IEEEtran}
\else
    \documentclass[11pt,draftcls,onecolumn]{IEEEtran}
\fi

\usepackage{amsthm}
\usepackage{amsmath}
\usepackage{fixmath}
\usepackage{amsfonts}
\usepackage[table]{xcolor}
\usepackage{adjustbox}
\ifCLASSOPTIONcompsoc
  \usepackage[nocompress]{cite}
\else
  \usepackage{cite}
\fi
\usepackage{graphicx}
\usepackage{algorithm}
\usepackage{algorithmic}
\usepackage{accents}
\usepackage{array,color}
\usepackage{subfig}
\usepackage{dsfont}
\usepackage{amsmath,bm}
\usepackage{supertabular}
\usepackage{enumitem}
\usepackage{upgreek}

\def\p{{\boldsymbol p}}
\def\x{{\boldsymbol x}}
\def\y{{\boldsymbol y}}
\def\b{{\mathbf b}}
\def\c{{\mathbf c}}
\def\xii{{\boldsymbol \xi}}
\def\h{{\mathbf h}}
\def\H{{\mathbf H}}
\def\A{{\mathbf A}}
\def\B{{\mathbf B}}

\def\pii{\bar{{\boldsymbol \pi}}}
\def\xx{{\underaccent{\bar}{\boldsymbol{x}}}}

\def\FF{{\underaccent{\bar}{\mathcal{F}}}}

\newcounter{proposition} 
\newenvironment{proposition}{ \vspace{0.1cm} \noindent {\bf Proposition
\theproposition.} \addtocounter{proposition}{-1}\refstepcounter{proposition}\em } { \addtocounter{proposition}{1}
\vspace{0.1cm} \normalfont }

\newcounter{remark} 
\newenvironment{remark}{\vspace{0.1cm} \noindent {\bf Remark \theremark.} \addtocounter{remark}{-1}\refstepcounter{remark}} {\addtocounter{remark}{1}
	\vspace{0.1cm} \normalfont }
    
\begin{document}

\title{Robust Approximate Dynamic Programming for Large-scale Unit Commitment with Energy Storages}

\author{Yu~Lan,\textit{~Member,~IEEE},
        Qiaozhu~Zhai,\textit{~Member,~IEEE},
        Xiaoming~Liu,\textit{~Member,~IEEE},
        and~Xiaohong~Guan,\textit{~Fellow,~IEEE}
\IEEEcompsocitemizethanks{
\IEEEcompsocthanksitem Y. Lan, Q. Zhai, X. Liu and X. Guan are with the Ministry of Education Key Lab for Intelligent Networks and Network Security, Xi'an Jiaotong University, Xi'an, China. 
E-mail: \{ylan, qzzhai, xmliu, xhguan\}@sei.xjtu.edu.cn.
}
}

\markboth{IEEE TRANSACTIONS ON POWER SYSTEMS}
{Lan \MakeLowercase{\textit{et al.}}: Robust Approximate Dynamic Programming for Large-scale Unit Commitment with Energy Storages}

\IEEEtitleabstractindextext{%
\begin{abstract}
The multistage robust unit commitment (UC) is of paramount importance for achieving reliable operations considering the uncertainty of renewable realizations. 
The typical affine decision rule method and the robust feasible region method may achieve uneconomic dispatches as the dispatch decisions just rely on the current-stage information.
Through approximating the future cost-to-go functions, the dual dynamic programming based methods have been shown adaptive to the multistage robust optimization problems, while suffering from high computational complexity.
Thus, we propose the robust approximate dynamic programming (RADP) method to promote the computational speed and the economic performance for large-scale robust UC problems. 
RADP initializes the candidate points for guaranteeing the feasibility of upper bounding the value functions, solves the linear McCormick relaxation based bilinear programming to obtain the worst cases, and combines the primal and dual updates for this hybrid binary and continuous decision-making problem to achieve fast convergence. 
We can verify that the RADP method enjoys a finite termination guarantee for the multistage robust optimization problems with achieving suboptimal solutions. 
Numerical tests on 118-bus and 2383-bus transmission systems have demonstrated that RADP can approach the suboptimal economic performance at significantly improved computational efficiency.
\end{abstract}

}

\maketitle

\IEEEdisplaynontitleabstractindextext

\IEEEpeerreviewmaketitle

\section*{Nomenclature}
\normalsize
\setlength{\tabcolsep}{0.8mm}
\begin{supertabular} {l p{0.75\linewidth}}
Variables        \\  
$x_{i,t}^{\rm{g}}$ & On/off status of the thermal unit $i$ at stage $t$ (0/1). \\
$x_{s,t}$ & Binary variable to represent that energy storage $s$ charge/discharge at stage $t$, $x_{s,t}=1$ denotes charging, $x_{s,t}=0$ denotes discharging. \\
$\boldsymbol{x}$ & Vector of all the binary variables, including $x_{i,t}^{\rm{g}}$ and $x_{s,t}$.\\
$p_{i,t}^{{\rm{g}}}$ & Power output of generator $i$ at stage $t$. \\   
$p_{s,t}^{\rm{sc}}$,$p_{s,t}^{\rm{sd}}$ & Charging/discharging power of the batteries. \\
$E_{s,t}$ &Storage level of battery $s$ at the end of $t$. \\
$\y_t$ & The $t$-stage dispatch decision vectors. \\
$p_{i,t}^{\rm{g,min}}, p_{i,t}^{\rm{g,max}}$ & The auxiliary variables to stand for the robust feasible regions of the unit dispatch decisions. \\
$E_{s,t}^{\rm{min}}, E_{s,t}^{\rm{max}}$ & The auxiliary variables to represent the robust feasible regions of the storage level decisions. \\
$p_{s,t}^{\rm{sc,min}}, p_{s,t}^{\rm{sc,max}}$ & The auxiliary variables for the robust feasible regions of the charging decisions. \\
$p_{s,t}^{\rm{sd,min}}, p_{s,t}^{\rm{sd,max}}$ & The auxiliary variables for the robust feasible regions of the discharging decisions. \\ [8pt]
Parameters  \\
$t$ & Index of stages, $t \in \mathcal{T}$. \\
$t_\delta$ & Time interval between stage $t$ and $t+1$.\\
$i$ & Index of the generators, $i \in \mathcal{N}_g$. \\
$s$ & Index of the batteries, $s \in \mathcal{N}_s$. \\
$r$ & Index of the renewable units, $r \in \mathcal{N}_r$. \\
$d$ & Index of the demand nodes, $d \in \mathcal{N}_d$.\\
$\ell$ & Index of the transmission lines, $\ell \in \mathcal{N}_\ell$.\\
$\mathcal{N}_g/\mathcal{N}_s/\mathcal{N}_r$ & Set of the generators/batteries/renewable units. \\
$\mathcal{N}_\ell/\mathcal{N}_d$ & Set of the transmission lines/demand nodes. \\
$n$ & Index of the state variable. \\
$N$ & Number of the dimensions for the state space. \\
$\kappa$ & Index of the supporting hyperplanes in outer approximations, $\kappa = 1,\cdots,\mathcal{K}$. \\
$C_{i}^{\rm{g}}$  & Fuel cost of the generators. \\
$C_{i}^{\rm{up}}$ & Start-up cost of the generators. \\
$\xi_{r,t}$ & Output of renewable unit $r$ at stage $t$. \\
$\xii_t$ & {Vector of renewable realizations at stage $t$.} \\
$\Upxi$ & Uncertainty set of the renewable outputs.\\
$p_{d,t}^{\rm{dm}}$ & Load demand at node $d$. \\
$\Gamma_{\ell,i}^{\rm g},\Gamma_{\ell,s}^{\rm{ss}}$ & Shift factor values of generators, storage units for transmission line $\ell$.\\
$\Gamma_{\ell,r}^{\rm r},\Gamma_{\ell,d}^{\rm{dm}}$ & Shift factor values of renewable units, loads for transmission line $\ell$.\\
$\bar{G}_\ell$ & The flow limit for transmission line $\ell$.\\
$\underaccent{\bar}{p}_i^{\rm{g}},\bar{p}_i^{\rm{g}}$ & Power output limits of generators.\\
${p}_i^{\rm{u}},p_i^{\rm d}$ & Startup/shutdown ramp limits of generators.\\
$\Delta_i^{\rm{u}},\Delta_i^{\rm{d}}$ & Ramp-rate up/down limits of generators. \\
$\underaccent{\bar}{E}_s, {{\bar E}_s}$ & Lower/upper bounds of battery storage levels.\\
${{\bar p}_s^{\rm{sc}}}$,${{\bar p}_s^{\rm{sd}}}$ & Maximum charging and discharging power of the batteries. \\
$\alpha^{\rm{sc}}, \alpha^{\rm{sd}}$ & Charging/discharging efficiency of batteries.\\
$\mathcal{X}$ & Feasible region of $\x$.\\
$\mathcal{Y}_t$ & The $t$-stage feasible region of $\y_t$. \\
$\underaccent{\bar}{y}_{n,t}, \bar{y}_{n,t}$ & Lower/upper limit for $t$-stage variable $y_{n,t}$. \\
$F^{\rm{uc}}$ & The UC-stage problem. \\
$F_{t}(\y_{t-1},\xii_{t})$ & The ED $t$-stage problem. \\
$\underaccent{\bar}{F}_{t}(\y_{t-1}, \xii_{t})$ & The $t$-stage lower bound problem.\\
$\bar{F}_t(\y_{t-1})$ & The $t$-stage upper bound problem.\\
$(\y_t^j,\bar{\mathcal{F}}_{t+1}(\y_t^j))$ & Candidate points used for convex combination, $j=1,\cdots, J_t$.\\
$\mathbb{Y}$ & The box to represent the feasible region of $\y_{t}$.\\
$\mathbb{S}$ & The Simplex with satisfying $\mathbb{Y} \subset \mathbb{S}$.\\
$\bar{F}_t^{\rm{bi}}(\pii_t,\xii_t)$ & The bilinear programming problem of upper bounding.\\
$\bar{F}_t^{\rm{McR}}(\pii_t, \xii_t)$ & The linear McCormick relaxation based upper bounding problem. \\
$\bar{F}_t^{\rm{pri}} (\y_{t-1},\xii_{t})$ & The primal upper bound with knowing $\xii_{t}$.\\
$\mathrm{env}(\cdot)$ & The lower convex envelope. \\
$\delta_{\underaccent{\bar}{\boldsymbol{y}}_{t}}(\y_{t})$ & The indicator function, $\delta_{\underaccent{\bar}{\boldsymbol{y}}_{t}}(\y_{t})=0$ if ${\y_{t}}=\underaccent{\bar}{\boldsymbol{y}}_{t}$ and $\delta_{\underaccent{\bar}{\boldsymbol{y}}_{t}}(\y_{t})=+\infty$ otherwise.\\
\end{supertabular}
\section{Introduction}
\label{intro}
\IEEEPARstart{R}{obust} unit commitment (UC) combining energy storage operation is of paramount significance for improving reliability against the uncertainty from the renewable realizations, 
which includes the commitment decision-making and the economic dispatch (ED) operation. 
The unit commitment decisions and operational adaptivity are highly related to the modeling ways of the multistage ED decisions \cite{lorca2016multistageTPS,bertsimas2012adaptive,lorca2016multistageRO}.      

Before the robust optimization (RO) was introduced to the UC problem, many typical stochastic programming (SP) methods have been proposed to manage the randomness of the renewable outputs, including the scenario trees and the chance-constrained modeling, see e.g., \cite{wang2008security,wu2007stochastic,zheng2014stochastic,wang2011chance,wu2014chance}, and references therein.  
Due to the added large number of scenarios and the probabilistic constraints, the SP methods suffer from high computational complexity for large-scale problems. 
To tackle the issues, RO has been proposed to solve the UC problem by relying on the uncertainty set while considering the feasibility under the infinite number of renewable realizations. The research on robust UC has become popular starting with the two-stage robust UC optimization \cite{bertsimas2012adaptive,jiang2011robust,zeng2013solving}. 
Assuming knowing the full knowledge of the future multistage realizations, the commitment decisions in the two-stage robust models can be made by guaranteeing the robustness.
However, this assumption is too strong to be satisfied in reality, as multistage ED decisions need to consider the so-termed nonanticipativity to achieve the sequential operations of the decision-making systems. 

To enforce the nonanticipativity of the ED decisions, a more accurate multistage robust formulation has been proposed for the UC problem, and then the multistage affine decision functions of uncertain parameters have been used to solve the large-scale problems for computational tractability \cite{lorca2016multistageRO}. 
The multistage affine decision rule method has many applications as in multi-period optimal power flow problem \cite{jabr2014robust}, real-time dispatch of the automatic generation control systems \cite{li2015adjustable}, decentralized cooperative operation in the distribution systems \cite{attarha2019affinely}, coordination of the electricity and the natural gas systems \cite{he2016robust}. 
Compared to restricting the ED decisions explicitly as the affine functions of the uncertainty, the implicit decision methods have been proposed to obtain the multistage robust feasible regions of the ED decisions through constructing the multistage nonanticipative constraints \cite{zhai2016transmission, cobos2018robust,cobos2018robustTSE,li2019multi,zhou2020multistage}.        
Both the multistage affine decision rule method and the robust feasible region method update the ED decisions based on the current-stage information, without considering the effects of the cost-to-go for the future stages. Thus, the ED decisions may result in uneconomic solutions for a multistage decision problem. 

Recently, the robust dual dynamic programming (RDDP) method has been proposed to solve the multistage RO problems iteratively through upper/lower bounding the worst-case cost-to-go functions  \cite{georghiou2019robust}.  
The idea of approximating the cost-to-go functions is from the stochastic dual dynamic programming (SDDP) for solving the multistage SP problems \cite{shapiro2011analysis,philpott2013solving,baucke2017deterministic}. 
The expected cost-go-to functions can also be evaluated by the approximate dynamic programming (ADP) method, which is to train the expected cost-go-to functions based on the preseclected structures using the suboptimal solution updates \cite{moazeni2018risk,shuai2018stochastic,xu2013kernel,zeng2018dynamic}.  
Compared to the probabilistic convergence of SDDP and the suboptimality of ADP, RDDP guarantees finite convergence to optimal solutions based on the deterministic updates of the candidate points. Many applications of RDDP/SDDP have been reported in multistage ED problems \cite{lu2019multi,papavasiliou2017application,shi2020enhancing,lan2022fast}, and multistage management problems in microgrids \cite{bhattacharya2016managing,shi2019multistage}.  
The extensions of SDDP/RDDP have been explored to solve the multistage stochastic/robust optimization with discrete recourse decisions as the stochastic dual dynamic integer programming (SDDiP) \cite{zou2019stochastic,zou2018multistage,hjelmeland2018nonconvex,ding2020multi} and the fast robust dual dynamic programming (FRDDP) \cite{xiong2022multi}.
However, the SDDiP method needs to binarize the continuous state variables to achieve valid, tight and finite cut generations, which is intractable for large-scale systems. The FRDDP method treats the binary decision variables as the continuous decision variables in the convex combinations of inner approximations, which is hard for upper bounding the value functions as in the systems with all the continuous recourse variables.
Therefore, it still needs to develop fast methods for solving the multistage RO problems with discrete and continuous decision variables.


The goal of the present work is to accelerate robust UC for large-scale systems by leveraging the latest advances in multistage robust dynamic programming techniques. Our contribution is two-fold. First, we develop the fast worst-case calculation based on the upper bounding of the worst-case cost-to-go function, which includes initializing the basic limited candidate points to guarantee the feasibility of convex combination and solving the bilinear programming under the linear McCormick relaxations.
Second, we further propose the robust approximate dynamic programming (RADP) scheme to solve the multistage robust UC, which can achieve suboptimal solutions with the finite termination guarantee. 
Specifically, we have successfully conducted the RADP for solving the robust UC problem with improved computational and economical performance based on both the primal and dual updates. The RADP enjoys efficient updates per iteration for both the commitment solutions and the worst-case cost-to-go approximations, with finite termination verified by our analysis and numerical tests. Numerical results also confirm the excellent performance of the RADP in achieving improved economical solutions by using the converged cost-to-go functions.

The rest of the paper is organized as follows. The multistage robust UC problem and the decoupled UC-stage and ED-stage formulations are introduced in Section \ref{section2}. Section \ref{section3} presents the RADP scheme along with the finite termination analysis for robust UC with hybrid binary and continuous decisions. Several numerical tests presented in Section \ref{section4} corroborate the faster computational performance of RADP relative to the RDDP solver, and the improved economical performance over the affine decision rule method and the robust feasible region method. The paper is wrapped up in Section \ref{section5}.

\section{Problem Formulation}
\label{section2}
Consider a multistage robust UC problem in the transmission systems as in \eqref{UC problem}, the objective is to minimize the commitment cost and the multistage worst-case dispatch cost:
\begin{subequations}\label{UC problem}
\begin{align}
    & \min_{\boldsymbol{x},\p(\cdot)} \sum_{i \in \mathcal{N}_g} C_i^{\mathrm{up}}(x_{i,\cdot}^{\rm{g}})+ \max_{\xii_{[t]} \in \Upxi} \sum_{t\in \mathcal{T}} \sum_{i \in \mathcal{N}_g} C_i^{\mathrm g}(p_{i,t}^{\mathrm g}(\xii_{[t]})) t_{\delta} \label{UC-obj} \\
    & \mathrm{s.t.} \;\; \x \in \mathcal{X}, \\
    & p_{i,t}^{\mathrm g}(\xii_{[t]})-p_{i,t-1}^{\mathrm g}(\xii_{[t-1]}) \leq \Delta_i^{\mathrm {u}} t_\delta x_{i,t-1}^{\rm{g}} + p_i^{\mathrm {u}}(x_{i,t}^{\rm{g}}-x_{i,t-1}^{\rm{g}}) \nonumber \label{UC-rampup}\\
    & \quad \quad \quad \quad +\bar{p}_i^{\mathrm g}(1-x_{i,t}^{\rm{g}}), \forall i\in \mathcal{N}_g, t\in \mathcal{T} \\
    & p_{i,t-1}^{\mathrm g}(\xii_{[t-1]})-p_{i,t}^{\rm g}(\xii_{[t]}) \leq \Delta_i^{\mathrm {d}} t_\delta x_{i,t}^{\rm{g}} + p_i^{\mathrm {d}}(x_{i,t-1}^{\rm{g}}-x_{i,t}^{\rm{g}}) \nonumber \label{UC-rampdown}\\
    & \quad \quad \quad \quad +\bar{p}_i^{\mathrm g}(1-x_{i,t-1}^{\rm{g}}), \forall i\in \mathcal{N}_g, t\in \mathcal{T} \\
    & E_{s,t}(\xii_{[t]})= E_{s,t-1}(\xii_{[t-1]})  + p_{s,t}^{\rm {sc}} (\xii_{[t]}) {\alpha^{\rm{sc}}} t_\delta \nonumber \\ & \quad \quad \quad \quad \quad \quad -p_{s,t}^{\rm {sd}}(\xii_{[t]}) t_\delta /{\alpha^{\rm{sd}}}, 
    \;s \in \mathcal{N}_s,\;t \in \mathcal{T} \label{UC-storage dynamics}\\
    & -\bar{G}_\ell \leq \sum_{i \in \mathcal{N}_g} \Gamma_{\ell, i}^{\rm g} p_{i,t}^{\rm g}(\xii_{[t]}) +\sum_{r \in \mathcal{N}_r} \Gamma_{\ell, r}^{\rm r} \xi_{r,t} \nonumber \\
    & +\sum_{s\in \mathcal{N}_s} \Gamma_{\ell, s}^{\rm {ss}} (p_{s,t}^{\rm{sd}}(\xii_{[t]})-p_{s,t}^{\rm{sc}}(\xii_{[t]})) -\sum_{n\in \mathcal{N}_d} \Gamma_{\ell,d}^{\rm{dm}} p_{d,t}^{\rm{dm}} \leq \bar{G}_\ell, \nonumber \\
    & \quad \quad \quad \quad \quad \quad \quad \quad \quad \quad \forall \ell \in \mathcal{N}_\ell, t\in \mathcal{T} \label{UC-lineflow}\\
    & \quad \quad \sum_{i \in \mathcal{N}_g} p_{i,t}^{\rm g}(\xii_{[t]}) +\sum_{s \in \mathcal{N}_s}(p_{s,t}^{\rm{sd}}(\xii_{[t]})-p_{s,t}^{\rm{sc}}(\xii_{[t]})) \nonumber \\
    & \quad \quad \quad + \sum_{r \in \mathcal{N}_r} \xi_{r,t}=\sum_{d \in \mathcal{N}_d} p_{d,t}^{\rm{dm}}, \forall t\in \mathcal{T} \label{UC-powerbalance}\\
    & \quad \quad x_{i,t}^{\rm{g}} \underaccent{\bar}{p}_i^{\rm g} \leq p_{i,t}^{\rm g}(\xii_{[t]}) \leq x_{i,t}^{\rm{g}} \bar{p}_i^{\rm g},\; \forall i \in \mathcal{N}_g, t \in \mathcal{T} \label{UC-unitlimits} \\
    & \quad \quad \quad \underaccent{\bar}{E}_s \le E_{s,t}(\xii_{[t]}) \le {{\bar E}_s},\;s \in \mathcal{N}_s,\;t \in \mathcal{T} \label{UC-storagelimits} \\ 
    & \quad \quad \quad \quad 0 \le p_{s,t}^{\rm{sc}}(\xii_{[t]}) \le \bar p_s^{\rm{sc}} x_{s,t},\;s \in \mathcal{N}_s,\;t \in \mathcal{T} \label{UC-chargelimits}\\
    & \quad \quad  0 \le p_{s,t}^{\rm{sd}}(\xii_{[t]}) \le \bar p_s^{\rm{sd}} (1-x_{s,t}),\;s \in \mathcal{N}_s,\;t \in \mathcal{T} \label{UC-dischargelimits} 
\end{align}
\end{subequations}
where $\x$ includes the binary variables $x_{i,t}^{\rm{g}}$ and $x_{s,t}$ satisfying $\x \in \mathcal{X}$. $x_{i,t}^{\rm{g}}$ stands for the commitment decisions determined before knowing the realizations of the uncertainty, which needs to satisfy all the commitment constraints as the start-up/shut-down constraints and the minimum up/down times constraints. 
The formulation details about $\x \in \mathcal{X}$ can be found in e.g. \cite{zhou2020multistage}. 
Here the commitment decisions are two-phase commitments, which means that $\x$ is not adaptive to the uncertainty realizations. 
The non-convex constraints to prevent simultaneous charging and discharging are represented by the binary variables $x_{s,t}$, which guarantee that the energy storages will not charge and discharge simultaneously.
The intertemporal constraints consist of the ramping up/down limits of the units in \eqref{UC-rampup}-\eqref{UC-rampdown} and the storage charging/discharging dynamics in \eqref{UC-storage dynamics}. 
The different charging and discharging coefficients are considered in \eqref{UC-storage dynamics}.   
The other constraints show transmission line flow limits in \eqref{UC-lineflow}, the power balance in \eqref{UC-powerbalance}, output limits of the units in \eqref{UC-unitlimits}, battery storage limits in \eqref{UC-storagelimits}, charging and discharging power limits in \eqref{UC-chargelimits} and \eqref{UC-dischargelimits}. To guarantee the sequential operation in reality, namely, the so-termed nonanticipativity, all the multistage ED decisions $\p_t(\xii_{[t]})$ are made based on the realizations of the uncertainty before stage $t$ as $\xii_{[t]}=(\xii_1,\cdots,\xii_t)$, without knowing the full knowledge of the future realizations.

The multistage RO problem in \eqref{UC problem} can be written into the following general form as:
\begin{align} \label{General problem}
\min & \;\c^\top\x + \max_{\xii_{[t]} \in \Upxi} \sum\nolimits_{t=1}^T \b_t^\top\y_t(\xii_{[t]})  \nonumber \\
\mathrm{s.t.} & \;\; \x \in \mathcal{X}, \nonumber \\
& \; \B_1\y_1(\xii_{1}) \geq \h_1(\x) + \H_1\xii_{1},\;\forall \; \xii_{1} \in \Upxi \nonumber \\
&\; \A_t\y_{t-1}(\xii_{[t-1]})+\B_t\y_t(\xii_{[t]}) \geq \h_t(\x)+\H_t\xii_{t},\nonumber \\
&\;\; \forall\; \xii_{t} \in \Upxi, \forall\; t=2,\cdots,T \nonumber \\
&\; \y_t(\xii_{[t]}) \in \mathcal{Y}_t,\xii_{t} \in \Upxi \; \mathrm{and} \; t=1,\cdots,T, 
\end{align}
where the vector $\y_t$ stands for the $t$-stage dispatch decision variables in \eqref{UC problem}, $\mathcal{Y}_t$ represents the $t$-stage feasible region of $\y_t$. $\x$ stands for the binary decisions in \eqref{UC problem}, which would be included in $\h_t$ for the constraints about $\y_t$. 

Assume that the uncertainty set is stagewise independent, the optimization problem \eqref{General problem} can be decoupled equivalently into $T+1$ two-stage subproblems as the UC-stage problem \eqref{stage 0 problem} and the $t$-stage ED problem \eqref{stage t problem}, $t=1,\cdots,T$:
\begin{align}
F^{\rm{uc}}=\min & \;\c^\top\x +\mathcal{F}_{1}(\y_0) \nonumber \\
\mathrm{s.t.} & \;\; \x \in \mathcal{X}, \label{stage 0 problem} \\
F_{t}(\y_{t-1},\xii_{t})&=\min \;\b_t^\top \y_t+\mathcal{F}_{t+1}(\y_t) \nonumber \\
\mathrm{s.t.} \;\A_t&\y_{t-1}+\B_t\y_t \geq \h_t(\x) + \H_t\xii_t, \nonumber \\
&\;\y_t \in \mathcal{Y}_t, \label{stage t problem}
\end{align}
where the UC-stage problem \eqref{stage 0 problem} denotes the commitment problem before the ED decision-making, the $t$-stage problem \eqref{stage t problem} denotes the ED problem in stage $t$ with given binary decisions $\x$, $\y_0$ is the initial state information.
Based on the $t$-stage problem $F_{t}(\y_{t-1},\xii_{t})$, the $t$-stage worst-case cost-to-go function $\mathcal{F}_{t}(\y_{t-1})$ is defined as: 
\begin{align}
    \mathcal{F}_{t}(\y_{t-1})=\max \{F_{t}(\y_{t-1},\xii_{t}) : \xii_{t} \in \Upxi_t\}.
\end{align}
One can calculate the optimal solutions of \eqref{stage t problem} with the information of the optimal worst-case cost-to-go functions $\mathcal{F}_{t+1}(\y_t)$. However, it is very challenging to obtain the optimal worst-case cost-to-go functions \cite{georghiou2019robust}. Thus, approximating the worst-case cost-to-go functions can help to solve the multistage RO problems.

Approximating the cost-to-go function through the outer approximation method as \cite{georghiou2019robust}, one can written the $t$-stage ED problem as: 
\begin{subequations} \label{stage t lower-bound problem}
\begin{align}%
\underaccent{\bar}{F}_{t}&(\y_{t-1}, \xii_{t})=\min \;\b_t^\top \y_t + \phi_{t+1}  \\
&\mathrm{s.t.} \;\A_t\y_{t-1}+\B_t\y_t \geq \h_t(\x) + \H_t\xii_t, \y_t \in \mathcal{Y}_t, \\
& \phi_{t+1} \geq \underaccent{\bar}{F}_{t+1}(\underaccent{\bar}{\y}_{t}^\kappa, \xii_{t+1}^\kappa) - \underaccent{\bar}{\boldsymbol{\pi}}_{t+1,\kappa}^{\top}\A_{t+1}(\y_{t}-\underaccent{\bar}{\y}_{t}^\kappa), \nonumber \\
&\; \quad \quad \quad \quad \quad \quad \quad \quad \quad \quad \quad \quad \quad \kappa =1, \cdots, \mathcal{K}, \label{lower-bound planes}
\end{align}
\end{subequations}
where \eqref{lower-bound planes} represents the supporting hyperplane based outer approximations for the cost-to-go functions, $\underaccent{\bar}{\boldsymbol{\pi}}_{t+1,\kappa}$ is the dual variable of solving $\underaccent{\bar}{F}_{t+1}({\y}_{t}, \xii_{t+1})$ when giving ${\y}_{t}=\underaccent{\bar}{\y}_{t}^\kappa $ and $\xii_{t+1}=\xii_{t+1}^\kappa$, $- \A_{t+1}^{\top}\underaccent{\bar}{\boldsymbol{\pi}}_{t+1,\kappa}$ denotes the $\kappa$-th subgradient of the functions. As the outer approximations are constructed by the supporting hyperplanes, the resultant values $\underaccent{\bar}{F}_{t} (\y_{t-1}, \xii_{t})$ in \eqref{stage t lower-bound problem} represent the lower bounds of the optimal worst-case cost-to-go functions. 

The so-termed RDDP scheme proposes to generate the worst cases $\xii_t$ used in \eqref{stage t lower-bound problem} according to the upper bounds of the worst-case cost-to-go functions \cite{georghiou2019robust}. 
Through refining the upper and lower bounds of the cost-to-go functions iteratively, RDDP obtains the optimal value functions when the upper and lower bounds converge. 
Due to the deterministic updates of the worst cases, RDDP guarantees finite convergence for achieving optimal solutions.  
Nevertheless, the iterative RDDP scheme suffers from high computational complexity for solving mixed-integer linear program (MILP) problems to obtain the worst cases, especially for large-scale optimization problems. Furthermore, RDDP solves the problems with all the continuous recourse decisions, without considering the RO problems with discrete decisions as \eqref{stage 0 problem}.  
Therefore, we propose the RADP method to solve the robust UC with discrete and continuous decision variables based on the primal and dual updates, which enjoys low computational complexity for large-scale systems.

\section{Robust Approximate Dynamic Programming}
\label{section3}
Thanks to the decoupled subproblems, one can update the worst-case cost-to-go functions based on the commitment solutions of the UC-stage problem. Given the binary decisions, the optimal worst-case cost-to-go functions about continuous ED decisions can be achieved by refining the lower and upper bounds. Nevertheless, the binary decisions are coupled in $T$ stages, which need to consider the influence of the sequential operations of the ED decisions. Thus, we first solve the UC-stage optimization problem by introducing robust nonanticipative constraints about the dispatch decisions. Furthermore, we propose to solve the relaxed upper-bounding problem with limited initialized candidate points for guaranteeing feasibility. Thirdly, we propose the RADP scheme combining both the primal and dual updates, and analyze the finite convergence for solving the multistage robust UC problem. 

\subsection{Nonanticipativity constrained UC-stage optimization}
To guarantee the robustness of the commitment decisions and the nonanticipativity of multistage ED decisions, we introduce the nonanticipative constraints about the generators and energy storages into the UC-stage problem as in the robust feasible region methods \cite{zhai2016transmission, cobos2018robust}:
\begin{align}
    & p_{i,t}^{\mathrm{g,max}}-p_{i,t-1}^{\mathrm {g,min}} \leq \Delta_i^{\mathrm u} t_\delta x_{i,t-1}^{\rm{g}} + p_i^{\mathrm u}(x_{i,t}^{\rm{g}}-x_{i,t-1}^{\rm{g}}) \label{nonanti-pg}\nonumber \\
    & \quad \quad \quad \quad \quad \quad +\bar{p}_i^{\mathrm g}(1-x_{i,t}^{\rm{g}}), \forall \; i\in \mathcal{N}_g, t\in \mathcal{T} \\
    & p_{i,t-1}^{\mathrm {g,max}}-p_{i,t}^{\rm {g,min}} \leq \Delta_i^{\mathrm d} t_\delta x_{i,t}^{\rm{g}} + p_i^{\mathrm d}(x_{i,t-1}^{\rm{g}}-x_{i,t}^{\rm{g}}) \nonumber \\
    & \quad \quad \quad \quad \quad \quad+\bar{p}_i^{\mathrm g}(1-x_{i,t-1}^{\rm{g}}), \forall \; i\in \mathcal{N}_g, t\in \mathcal{T} \\
    & x_{i,t}^{\rm{g}} \underaccent{\bar}{p}_i^{\rm g} \leq p_{i,t}^{\rm{g,min}}\leq p_{i,t}^{\rm g}(\xii_{[t]}) \leq p_{i,t}^{\rm{g,max}} \leq x_{i,t}^{\rm{g}} \bar{p}_i^{\rm g},\nonumber \\ 
    & \quad \quad \quad \quad \quad \quad \quad \quad\quad \quad \quad \forall \; i \in \mathcal{N}_g, t \in \mathcal{T} \\
    & 0 \le p_{s,t}^{\rm{sc,min}}\le p_{s,t}^{\rm{sc}}(\xii_{[t]}) \le p_{s,t}^{\rm{sc,max}} \le \bar p_s^{\rm{sc}} x_{s,t}, \nonumber \\ 
    & \quad \quad \quad \quad \quad \quad \quad \quad \quad \quad \quad \forall\;s \in \mathcal{N}_s,\;t \in \mathcal{T} \\
    & 0 \le p_{s,t}^{\rm{sd,min}}\le p_{s,t}^{\rm{sd}}(\xii_{[t]}) \le p_{s,t}^{\rm{sd,max}}\le \bar p_s^{\rm{sd}} (1-x_{s,t}), \nonumber \\ 
    & \quad \quad \quad \quad \quad \quad \quad \quad \quad \quad \quad \forall \;s \in \mathcal{N}_s,\;t \in \mathcal{T} \\
    & E_{s,t}^{\rm{min}}= E_{s,t-1}^{\rm{min}}  + p_{s,t}^{\rm {sc,min}} {\alpha^{\rm{sc}}} t_\delta -p_{s,t}^{\rm {sd,max}} t_\delta /{\alpha^{\rm{sd}}}, \nonumber \\ 
    & \quad \quad \quad \quad \quad \quad \quad \quad \quad \quad \quad \forall \;s \in \mathcal{N}_s,\;t \in \mathcal{T} 
\end{align}
\begin{align}
    & E_{s,t-1}^{\rm{max}}= E_{s,t}^{\rm{max}} + p_{s,t}^{\rm {sc,max}} {\alpha^{\rm{sc}}} t_\delta -p_{s,t}^{\rm {sd,min}} t_\delta /{\alpha^{\rm{sd}}}, \nonumber \\ 
    & \quad \quad \quad \quad \quad \quad \quad \quad \quad \quad \quad \forall \;s \in \mathcal{N}_s,\;t \in \mathcal{T} \\
    & \quad \underaccent{\bar}{E}_s \le E_{s,t}^{\rm{min}} \le E_{s,t}^{\rm{max}} \le {{\bar E}_s},\forall\;s \in \mathcal{N}_s,\;t \in \mathcal{T}, \label{nonanti-et}
\end{align}
where $p_{i,t}^{\rm{g, min}}, p_{i,t}^{\rm{g, max}}, E_{s,t}^{\rm{min}}, E_{s,t}^{\rm{max}}, p_{s,t}^{\rm{sc, min}}, p_{s,t}^{\rm{sc, max}}, p_{s,t}^{\rm{sd, min}}, \\ p_{s,t}^{\rm{sd, max}}$ are auxiliary variables to stand for the robust feasible regions of the multistage dispatch decisions. 
They can be determined in the UC-stage optimization and directly passed to the ED stages for approximations of the cost-to-go functions. 
The robust feasible region methods try to obtain feasible region limits for the dispatch variables, which can be robust to any realization of the uncertainty. At the same time, the nonanticipative constraints help to satisfy the sequential dispatch property, which means the ED decisions can be made based on the current-stage realizations not on future realizations.

According to the nonanticipative constraints added for the multistage dispatch decisions, the UC-stage problem can be rewritten into \eqref{nonancipative UC-stage} as:
\begin{align} \label{nonancipative UC-stage}
F^{\rm{uc}}=&\min \;\c^\top\x + \phi_0 \nonumber \\
\mathrm{s.t.} & \; \x \in \mathcal{X},  \nonumber\\
& \phi_0 \geq \sum\nolimits_{t=1}^T \b_t^\top\y_t(\xii_{t}), \forall \; \xii_t \in \mathcal{U}_t, \nonumber\\
& \eqref{UC-lineflow}-\eqref{UC-powerbalance},\; \eqref{nonanti-pg}-\eqref{nonanti-et}, \forall \; \xii_t \in \mathcal{U}_t,
\end{align}
where $\mathcal{U}_t$ stands for a subset of $\Upxi_t$, which can be expanded by representative scenarios. The worst-case scenarios added into $\mathcal{U}_t$ can be calculated through solving the optimization problems in the ED stages. And the associated generations of the recourse decision variables in \eqref{nonancipative UC-stage} are based on the column-and-constraint generation (CCG) algorithm \cite{zeng2013solving}. 
Knowing the solutions $\x$ of \eqref{nonancipative UC-stage}, one can approximate the worst-case cost-to-go functions of the continuous ED decisions through upper/lower bounding. 


\subsection{Generations of worst cases through upper bounding}
\label{initial setting of the algorithm} 
Different from lower bounding the value functions through the outer approximations as \eqref{stage t lower-bound problem}, the upper bounds of the worst-case cost-to-go functions are constructed by the convex combination based inner approximation methods, which can be formulated as \cite{georghiou2019robust}:
\begin{subequations}\label{upper primal problem}
\begin{align}
    \bar{F}_t(\y_{t-1})&= \max_{\xii_{t} \in \Upxi_t} \min_{\boldsymbol{y}_t, \boldsymbol{\lambda}} \;\b_t^\top \y_t+\sum\nolimits_j {\lambda}_j \bar{\mathcal{F}}_{t+1}(\y_t^j)  \label{upper obj} \\
   \rm{s.t.}\; &  \B_t \y_t \geq \h_t(\x)-\A_t\y_{t-1}+\H_t\xii_t, \;\y_t \in \mathcal{Y}_t \label{upper constraint b} \\
   & \boldsymbol{y}_t =\;\sum\nolimits_j {\lambda}_j \boldsymbol{y}_t^j,\; {\lambda}_j\geq 0,\; \sum\nolimits_j {\lambda}_j=1, \label{upper constraint c}
\end{align}
\end{subequations}
where the decision variables and the cost-to-go values are formed by the convex combinations of the candidate points for the upper-bound functions as in \eqref{upper obj} and \eqref{upper constraint c}, the candidate points $\{(\y_t^j,\bar{\mathcal{F}}_{t+1}(\y_t^j)), j=1,\cdots, J_t\}$ are calculated from the previous iterations. 
The worst cases can be generated by solving the upper-bound problems in \eqref{upper primal problem}.

However, the convex combination cannot guarantee the feasibility in the state space when the previously generated candidate points are limited. For example, $\y_t=\lambda_1 \y_t^1$ if $J_t=1$. In this case, the solution $\y_t$ can only be $\y_t^1$ based on \eqref{upper constraint c}. It may be infeasible for constraint \eqref{upper constraint b} as the newly updated $\y_{t-1}$ can be different from the previous one. The basic way to guarantee the feasibility of the convex combination is to choose all the extreme points of the state space as the candidate points \cite{philpott2013solving}. Assuming the dimension of the state is $N$, it will add $2^N$ variables for \eqref{upper constraint c} as the $2^N$ extreme points. It suffers from computational difficulty when $N$ is large. 
Recently, slack variables $\boldsymbol{v}_t^1 $ and $\boldsymbol{v}_t^2$ are introduced to describe $\boldsymbol{y}_t$ as $\boldsymbol{y}_t =\;\sum\nolimits_j {\lambda}_j \boldsymbol{y}_t^j + \boldsymbol{v}_t^1-\boldsymbol{v}_t^2$ to guarantee the feasibility using the previously generated points, with adding a new term $\boldsymbol{\tau}_t^\top (\boldsymbol{v}_t^1+\boldsymbol{v}_t^2)$ in the objective \cite{shi2020enhancing}. 
Nevertheless, it needs to tune the penalty coefficients $\boldsymbol{\tau}_t$ according to the subgradients of the lower-bound objectives. This is hard for guaranteeing the lower/upper bounding property as $\FF_t(\y_{t-1}) \leq \mathcal{F}_t(\y_{t-1}) \leq \bar{\mathcal{F}}_t(\y_{t-1})$, which is the key point for convergence to the optimal solutions \cite{georghiou2019robust}. 


Therefore, we propose to initialize the candidate points for upper bounding using the limited $N+1$ points as \eqref{initial points} when the dimension of the state is $N$:
\begin{align}\label{initial points}
    \{\beta\mathbf{e}_1 ,\cdots, \beta\mathbf{e}_n,\cdots, \beta \mathbf{e}_N, \boldsymbol{0}\},
\end{align}
where $\{\mathbf{e}_n\}_{n=1}^N$ denote the canonical basis of $\mathbb{R}^N$, $\beta=\sum_{n=1}^N \bar{y}_{n,t}$ stands for the summation of the states' upper limits. 
\begin{figure}[!t] 
\centering 
\includegraphics[trim = 290mm 127mm 80mm 80mm, clip,width=0.5\textwidth]{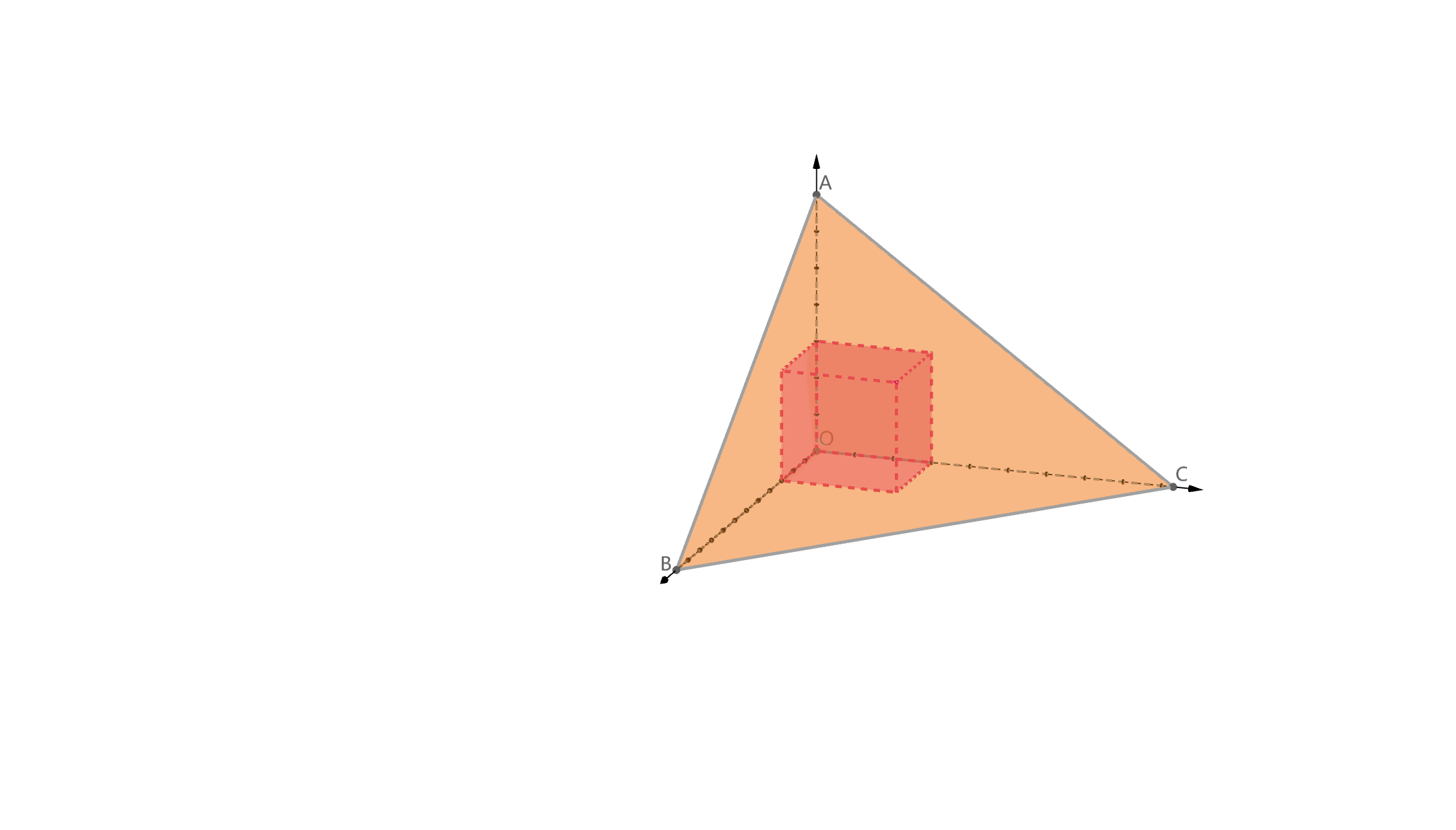} 
\caption{{Initialized $N+1$ candidate points for upper-bounding.}} 
\label{fig.simplex} 
\end{figure}

\begin{proposition}\label{upperbound-feasible}
{Assume the dimension of the state space is $N$, the initialized $N+1$ points in \eqref{initial points} can guarantee the feasibility of the convex combination for upper bounding in \eqref{upper primal problem}.} 
\end{proposition}
\begin{proof}
{The feasible region of $\y_t$ can be represented by the $N$-dimensional box $\mathbb{Y}$ as $\{0 \leq y_{n,t} \leq \bar{y}_{n,t},\; n=1,\cdots,N \} $. One can easily construct a Simplex $\mathbb{S}$ with $\{y_{n,t} \geq 0, \;  \sum_{n=1}^N y_{n,t} \leq \sum_{n=1}^N \bar{y}_{n,t}$\}, which satisfies $\mathbb{Y} \subset \mathbb{S}$. Fig. \ref{fig.simplex} shows the relationship between $\mathbb{Y}$ (the pink box) and $\mathbb{S}$ (the OABC) in 3-dimensional space. Finding any feasible points in $\mathbb{Y}$ can be achieved by finding any feasible points in $\mathbb{S}$ as $\mathbb{Y} \subset \mathbb{S}$. According to the convex analysis, any points in $\mathbb{S}$ can be represented by convex combination using the $N+1$ vertices of $\mathbb{S}$ as listed in \eqref{initial points}. One can complete the proof in a similar way when the variables have lower limits as $y_{n,t} \geq \underaccent{\bar}{y}_{n,t}$.  }
\end{proof}

Compared with using $2^N$ points to guarantee feasibility, the initialized limited $N+1$ candidate points result in much fewer variables for convex combinations, which is much more computationally efficient for upper bounding. 
Compared with introducing slack variables to guarantee the feasibility and tuning the penalized coefficients in the objectives, the bounding property $\FF_t(\y_{t-1}) \leq \mathcal{F}_t(\y_{t-1}) \leq \bar{\mathcal{F}}_t(\y_{t-1})$ is always guaranteed by the proposed initial $N+1$ candidate points.
Note that the $N+1$ points in \eqref{initial points} are the initialized candidate points for upper-bounding. The set of candidate points will be expanded based on the update of cost-to-go functions. 
With guaranteeing the feasibility of the inner approximations, the max-min upper-bound problem in \eqref{upper primal problem} can be solved to obtain the worst cases. 



Dualizing the inner minimization problem, the upper-bound problem \eqref{upper primal problem} becomes a maximization problem with a bilinear term in the objective as: 
\begin{align}\label{original bilinear programming}
    \bar{F}_t^{\rm{bi}}(\pii_t,\xii_t) &=  \max_{\xii_t,\pii_t} \; (\bar{\h}_t(\x)-\bar{\A}_t\y_{t-1})^\top \pii_t +\boldsymbol{\eta}_t^\top \xii_t \nonumber \\
     \mathrm{s.t.} \; \;\; &  \bar{\B}_t^\top \pii_t \leq \bar{\b}_t ,\; \boldsymbol{\eta}_t=\bar{\H}_t^\top \pii_t, \pii_t \geq 0, \xii_{t} \in \Upxi_t, 
\end{align}
where $\pii_t$ is the dual variables for the upper-bound problem \eqref{upper primal problem}, $\bar{\cdot}$ denotes the matrices/vectors for the upper-bound problem considering the inner approximations. Introducing auxiliary binary variables and adopting a Big-M reformulation, this maximization problem can be transformed into a MILP problem \cite{zeng2013solving,lee2013modeling}. The MILP-based methods are powerful to obtain the optimal solutions for \eqref{upper primal problem}. 
Nevertheless, the MILPs are time-consuming for large-scale system optimization problems, especially in an iterative scheme for achieving the optimal cost-to-go functions. 

The bilinear term $\boldsymbol{\eta}_t^\top \xii_t$ in \eqref{original bilinear programming} can be rewritten into $\sum_{r=1}^{\mathcal{N}_r} \eta_{r,t} \xi_{r,t}$. The $r$-th term $\eta_{r,t}\xi_{r,t}$ can be bounded by the convex and concave envelopes as \eqref{McEnve}, the so-termed linear McCormick relaxations \cite{mccormick1976computability,ben2017convex,deng2021optimal}, where $\underaccent{\bar}{\xi}_{r,t},\bar{\xi}_{r,t}, \underaccent{\bar}{\eta}_{r,t}, \bar{\eta}_{r,t}$ are the lower/upper limits for $\xi_{r,t}$ and $\eta_{r,t}$.
\begin{subequations}\label{McEnve}
\begin{align}
    & \theta_{r,t} \geq \eta_{r,t} \underaccent{\bar}{\xi}_{r,t} + \underaccent{\bar}{\eta}_{r,t} \xi_{r,t} - \underaccent{\bar}{\eta}_{r,t} \underaccent{\bar}{\xi}_{r,t}, \label{McEnve-a}\\
    & \theta_{r,t} \geq \eta_{r,t} \bar{\xi}_{r,t} + \bar{\eta}_{r,t} \xi_{r,t} - \bar{\eta}_{r,t} \bar{\xi}_{r,t}, \\
    & \theta_{r,t} \leq \eta_{r,t} \underaccent{\bar}{\xi}_{r,t} + \bar{\eta}_{r,t} \xi_{r,t} - \bar{\eta}_{r,t} \underaccent{\bar}{\xi}_{r,t}, \\
    & \theta_{r,t} \leq \eta_{r,t} \bar{\xi}_{r,t} + \underaccent{\bar}{\eta}_{r,t} \xi_{r,t} - \underaccent{\bar}{\eta}_{r,t} \bar{\xi}_{r,t}. \label{McEnve-d} 
\end{align}
\end{subequations}

Due to the linear relaxations, the upper-bound bilinear programming turns out to be an upper-bound linear programming (LP) problem $\bar{F}_t^{\rm{McR}}$ as 
\begin{align}\label{LP-based bilinear program}
    \bar{F}_t^{\rm{McR}}  (\pii_t,& \xii_t) = \max_{\xii_t,\pii_t}\; (\h_t(\x)-\A_t\y_{t-1})^\top \pii_t + \sum_{r=1}^{\mathcal{N}_r} \theta_{r,t} \nonumber \\
    \mathrm{s.t.} \; 
    &\; \B_t^\top \pii_t \leq \hat{\b}_t ,\; \boldsymbol{\eta}_t=\H_t^\top \pii_t, \pii_t \geq 0, \xii_{t} \in \Upxi_t \nonumber \\
    &\; \eqref{McEnve-a} \sim \eqref{McEnve-d},\; r \in \mathcal{N}_r.
\end{align}
Solving $\bar{F}_t^{\rm{McR}}$ can help to obtain the suboptimal solution for bilinear programming, which is faster than solving MILPs. 
Based on the commitment solutions of $F^{\mathrm{uc}}$ in \eqref{nonancipative UC-stage} and the worst-case solutions of $\bar{F}_t^{\rm{McR}}$, we propose the RADP scheme to solve the multistage robust UC problem. Furthermore, we will show that RADP will terminate in finite time achieving the converged suboptimal worst-case cost-to-go functions.   

\subsection{The RADP scheme}
\label{LP-based RDDP}
The RADP method is tabulated in Algorithm \ref{alg:RADP_UC}. In each iteration, the nonanticipative constraints based UC-stage problem is first to be solved to obtain the commitment decisions $\xx$, which are passed to the inner loop for calculating the representative scenarios and approximating the worst-case cost-to-go functions.

Initializing the upper/lower bounds as $\bar{\mathcal{F}}_t^0(\y_{t-1})$ and $\underaccent{\bar}{\mathcal{F}}_t^0(\y_{t-1})$, the inner loop includes the forward pass for generating the candidate points and the backward pass for refining the upper/lower bounds of the value functions.
The initialized upper and lower bounds can be set directly as $+\infty$ and $-\infty$. Or the bounds can be calculated by solving LP problems while viewing $\y_{t-1}$ as variables \cite{lan2022fast}. 
Although all the worst cases are generated by $\bar{F}_t^{\rm{McR}}$, the $\xii_t^{\rm{bw}}$ in the backward pass may be different from the $\xii_t^{\rm{fw}}$ in the forward pass as the upper-bound approximations in the backward pass have been updated. 
The basic idea of updating the bounds is to narrow the gaps between the lower and upper bounds until they converge. 
Note that the upper bounds are not refined using the values of $\bar{F}_t^{\rm{McR}}(\pii_t,\xii_t^{\rm{bw}})$ as they are just the approximate values based on relaxations. 
With the information $\xii_t^{\rm{bw}}$, the upper-bound value at $\underaccent{\bar}{\boldsymbol{y}}_{t-1}$ is updated by calculating $\bar{F}_t^{\rm{pri}} (\underaccent{\bar}{\y}_{t-1},\xii_t^{\rm{bw}})$ in \eqref{upper value problem}, which is the valid upper bound compared to $\bar{F}_t^{\rm{McR}}(\pii_t,\xii_t^{\rm{bw}})$.
When the inner loop converges with achieving $\underaccent{\bar}{\mathcal{F}}_1({\y}_0)=\bar{\mathcal{F}}_1(\y_0)$, all the worst cases $\xii_t^{\rm{bw}}$ are collected and added into $\mathcal{U}_t$ for the UC-stage problem in the next iteration to update $\xx$. 
\begin{align} \label{upper value problem}
    \bar{F}_t^{\rm{pri}} & (\y_{t-1},\xii_{t})= \min_{\boldsymbol{y}_t} \;\b_t^\top \y_t+\sum\nolimits_j {\lambda}_j \bar{\mathcal{F}}_{t+1}(\y_t^j)  \nonumber \\
   \rm{s.t.}\; &  \B_t \y_t \geq \h_t(\x)-\A_t\y_{t-1}+\H_t\xii_t, \; \eqref{upper constraint c}, \;\y_t \in \mathcal{Y}_t 
\end{align}

The RADP stops when the commitment solutions $\xx$ are the same in consecutive iterations, which means the multistage worst-case cost-to-go functions have converged in the inner loop, and $\mathcal{U}_t$ will not update anymore. 

\begin{algorithm}[!t]
\caption{RADP Scheme}
\label{alg:RADP_UC}
\begin{algorithmic}[1]
\REQUIRE Maximum iteration number $\mathcal{M}$, time horizon $T$.
\ENSURE  $\xx$, $\FF_t(\y_{t-1})$ and $\bar{\mathcal{F}}_{t}(\y_{t-1})$ for $t=1,\cdots,T$
\FOR{$m_1=1$ to $\mathcal{M}$}
\STATE Solve the UC-stage problem \eqref{nonancipative UC-stage} to obtain $\xx^{(m_1)}$. \label{Stage0 line}\\
If $\xx^{(m_1)}=\xx^{(m_1-1)}$, terminate, otherwise, go on. 
\STATE   Set $\underaccent{\bar}{\mathcal{F}}_t(\y_{t-1})=\underaccent{\bar}{\mathcal{F}}_t^0(\y_{t-1})$, $\bar{\mathcal{F}}_t(\y_{t-1})=\bar{\mathcal{F}}_t^0(\y_{t-1})$, $\underaccent{\bar}{\mathcal{F}}_{T+1}(\y_T)=\bar{\mathcal{F}}_{T+1}(\y_T)=0$.
\FOR{$m_2=1$ to $\mathcal{M}$} 
\STATE Forward pass:\\
For $t=1,\cdots,T$, let $\xii_t^{\rm{fw}}$ be optimal solutions of $\bar{F}_t^{\rm{McR}}(\pii_t,\xii_t)$ in \eqref{LP-based bilinear program}, obtain optimal solutions $\underaccent{\bar}{\y}_t$ for solving $\underaccent{\bar}{F}_t(\underaccent{\bar}{\y}_{t-1};\xii_t^{\rm{fw}})$ in \eqref{stage t lower-bound problem}.
\STATE Backward pass: \\ 
For stage $t=T,\cdots,1$, let $\xii_t^{\rm{bw}}$ be optimal solutions to $\bar{F}_t^{\rm{McR}}(\pii_t,\xii_t)$ in \eqref{LP-based bilinear program}.\\ 
Solving the primal upper-bound problems $\bar{F}_t^{\rm{pri}} (\underaccent{\bar}{\y}_{t-1},\xii_t^{\rm{bw}})$ in \eqref{upper value problem}, update the upper bounds as $\bar{\mathcal{F}}_t(\y_{t-1}) \leftarrow \mathrm{env} \big (\min \{\bar{\mathcal{F}}_t(\y_{t-1}), \bar{{F}}^{\rm{pri}}_t(\underaccent{\bar}{\boldsymbol{y}}_{t-1},\xii_t^{\rm{bw}})+\delta_{\underaccent{\bar}{\boldsymbol{y}}_{t-1}} (\y_{t-1})\}\big )$.  \\
Solving the dual problems of $\underaccent{\bar}{F}_t(\underaccent{\bar}{\y}_{t-1},\xii_t^{\rm{bw}})$ in \eqref{stage t lower-bound problem} to obtain $\underaccent{\bar}{\boldsymbol{\pi}}_t$. Update the lower bounds as \eqref{lower-bound planes}. \\
If $\underaccent{\bar}{\mathcal{F}}_1({\y}_0)=\bar{\mathcal{F}}_1(\y_0)$, update the set $\mathcal{U}_t=\mathcal{U}_t\cup \xii_t^{\rm{bw}}$, and go to Step \ref{Stage0 line}. 
\ENDFOR
\ENDFOR 
\RETURN $\xx$, $\FF_t(\y_{t-1})$ and $\bar{\mathcal{F}}_{t}(\y_{t-1})$ for $t=1,\cdots,T$
\end{algorithmic}
\end{algorithm}

\begin{remark}\label{convergence}
The RADP scheme terminates in finite time.
\end{remark}
\begin{proof}
The finite termination of the RADP scheme can be proved in two steps. First, the inner loop for the cost-to-go approximations can converge in finite iterations for given $\xx$. Second, the UC-stage problem generates a finite number of $\xx$. 

The sketch of the proof for the finite termination of the inner loop follows from \cite{georghiou2019robust}. After showing that there will have finite lower and upper bounds generated by the RADP, one can prove at least one of the bounds must be refined in the backward pass. 

Due to the convex piecewise-affine outer approximations, $\underaccent{\bar}{{F}}_t(\y_{t-1},\xii_{t})$ are linear programs with finite solutions, which generate finite supporting hyperplanes for $\FF_t$ under each $\xii_t$. Thus, there are finite lower bounds to be refined. With finite approximations $\FF_{t+1}$ and finite scenarios $\xii_t^{\rm{fw}}$ or $\xii_t^{\rm{bw}}$ generated by $\bar{F}_t^{\rm{McR}}(\cdot)$, finite candidate points $\underaccent{\bar}{\boldsymbol {y}}_t$ can be obtained from $\underaccent\bar{F}_t(\underaccent{\bar}{\boldsymbol {y}}_{t-1},\xii_t)$. Therefore, there are finite upper bound problems to solve with finite candidate points $\underaccent{\bar}{\boldsymbol {y}}_{t-1}$, which result in finite upper bounds to be refined.

Following the proof of contradiction by assuming none of the bounds were refined, which will imply that the inner loop will terminate at the beginning without refining any bounds.
By backward induction, for stage $t \leq T$, one can easily obtain\\
$\FF_t(\underaccent{\bar}{\boldsymbol y}_{t-1}) =\underaccent\bar{F}_t(\underaccent{\bar}{\boldsymbol y}_{t-1},\xii_t^{\rm{bw}})
\leq \bar{{F}}_t^{\rm{pri}}(\underaccent{\bar}{\boldsymbol y}_{t-1},\xii_t^{\rm{bw}})
=\bar{\mathcal{F}}_t(\underaccent{\bar}{\boldsymbol y}_{t-1}),$
where the equalities hold for the bounds are not refined, the inequality would be strengthened to be equality because none of the bounds need to be refined as $\FF_{t+1}(\y_{t}) = \mathcal{F}_{t+1}(\y_{t}) = \bar{\mathcal{F}}_{t+1}(\y_{t})$. It yields $\FF_{t}(\y_{t-1}) = \bar{\mathcal{F}}_{t}(\y_{t-1})$, which means the inner scheme terminates at the beginning. Therefore, at least one of the finite bounds must be refined in the backward pass.

Benefiting from the finite termination of the inner scheme, finite scenarios $\xii_t^{\rm{bw}}$ will be generated and added to the UC-stage optimization problem. Thus, the UC-stage problem has a finite number of solutions $\xx$ in the iterative scheme. 
\end{proof}

Due to the finite updates of the cost-to-go functions in the inner loop and the finite scenarios added to the UC-stage problem, the RADP scheme can achieve finite termination. 
RADP combines the primal scenario-based recourse decision variable updates for the UC-stage optimization and the dual updates for approximating the cost-to-go functions. 
Apart from the commitment solutions, one can obtain the converged worst-case cost-to-go functions when the scheme terminates. Therefore, given the specific realizations, the ED decisions can be calculated by solving $\underaccent{\bar}{F}_{t}(\y_{t-1}, \xii_{t})$ in \eqref{stage t lower-bound problem} or solving $\bar{F}_t^{\rm{pri}} (\y_{t-1},\xii_{t})$ in \eqref{upper value problem} as the lower/upper bounds converge to $\mathcal{F}_{t+1}(\y_{t})$. 
Our numerical tests have shown that the RADP method will achieve economical improvement by using the converged cost-to-go functions.

The RADP method is an approximate method for solving the large-scale robust UC problem, the reasons are two-fold: 1) the nonanticipative constraints used in the UC-stage problem \eqref{nonancipative UC-stage} are sufficient conditions for obtaining robust solutions, 2) the generated worst cases are approximate solutions based on the linear McCormick relaxations.  
Note that the termination analysis of the RADP scheme is based on the relatively complete recourse assumptions, which can be satisfied by introducing auxiliary variables in the constraints and adding the penalization terms in the objective functions \cite{zou2018multistage,zou2019stochastic}.

\section{Numerical Results}
\label{section4}
The proposed RADP and the comparison methods have been tested on a Windows server with Intel\textsuperscript{\textregistered} 8-core CPU @ 3.6 GHz (64GB RAM) in the MATLAB\textsuperscript{\textregistered} R2021b simulator. All the optimization problems have been modeled by Yalmip and solved by Gurobi v9.5.2. The 118-bus and 2383-bus test cases are implemented in the MATPOWER toolbox \cite{zimmerman2010matpower}. The basic information of the test cases are listed in Table \ref{table:characteristics of cases}, the power limits and other parameters are from \cite{zimmerman2010matpower,zhou2020multistage,li2019multi}. The dimension of the state space for each test case is the summation of the number of units and the number of energy storages for the intertemporal constraints in \eqref{UC problem}. 

\begin{figure*}[th]
\centering
    \subfloat{\label{Nr10_fig}\includegraphics[width=0.32\textwidth]{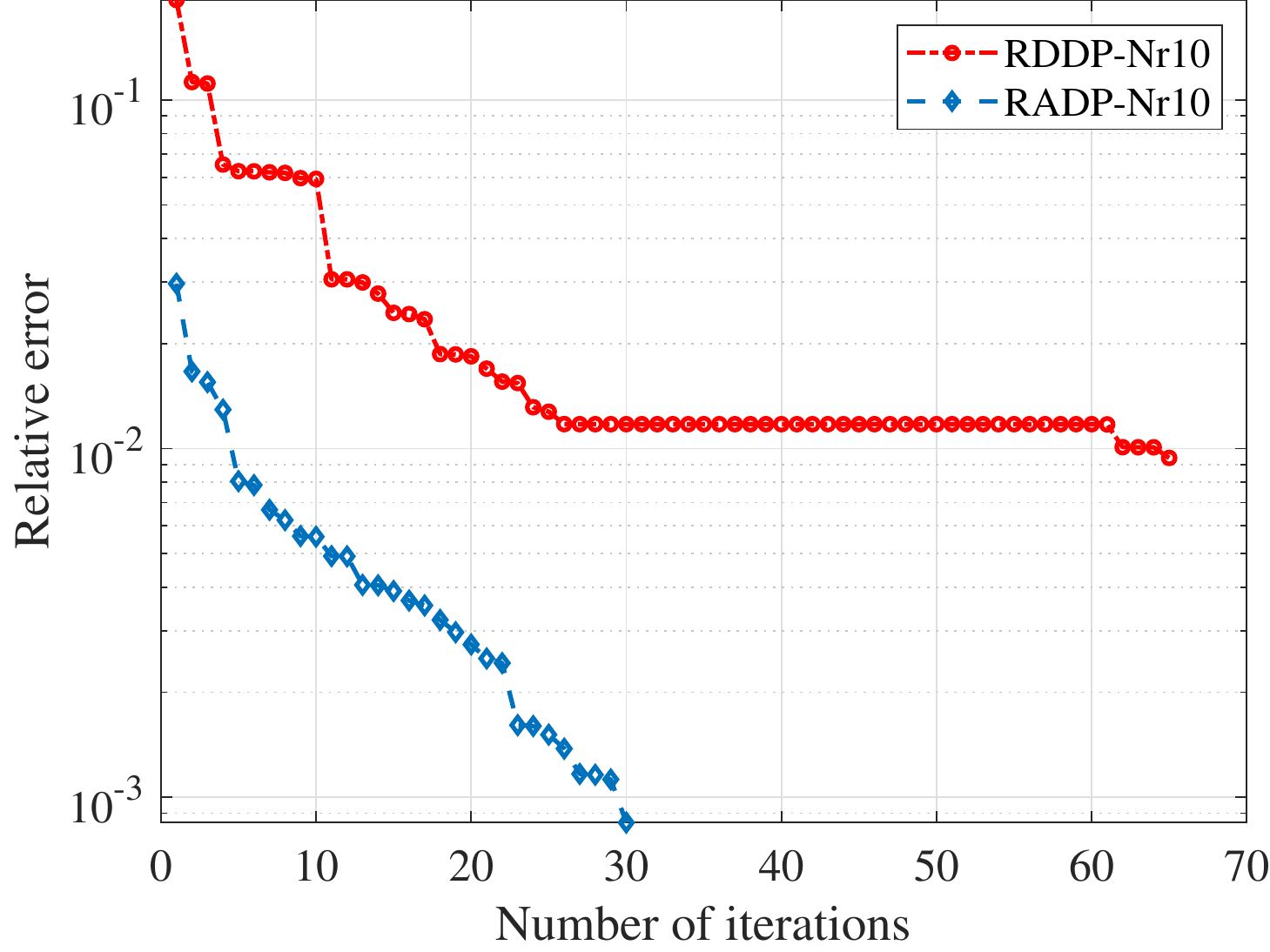}} 
    \subfloat{\label{Nr20_fig}\includegraphics[width=0.32\textwidth]{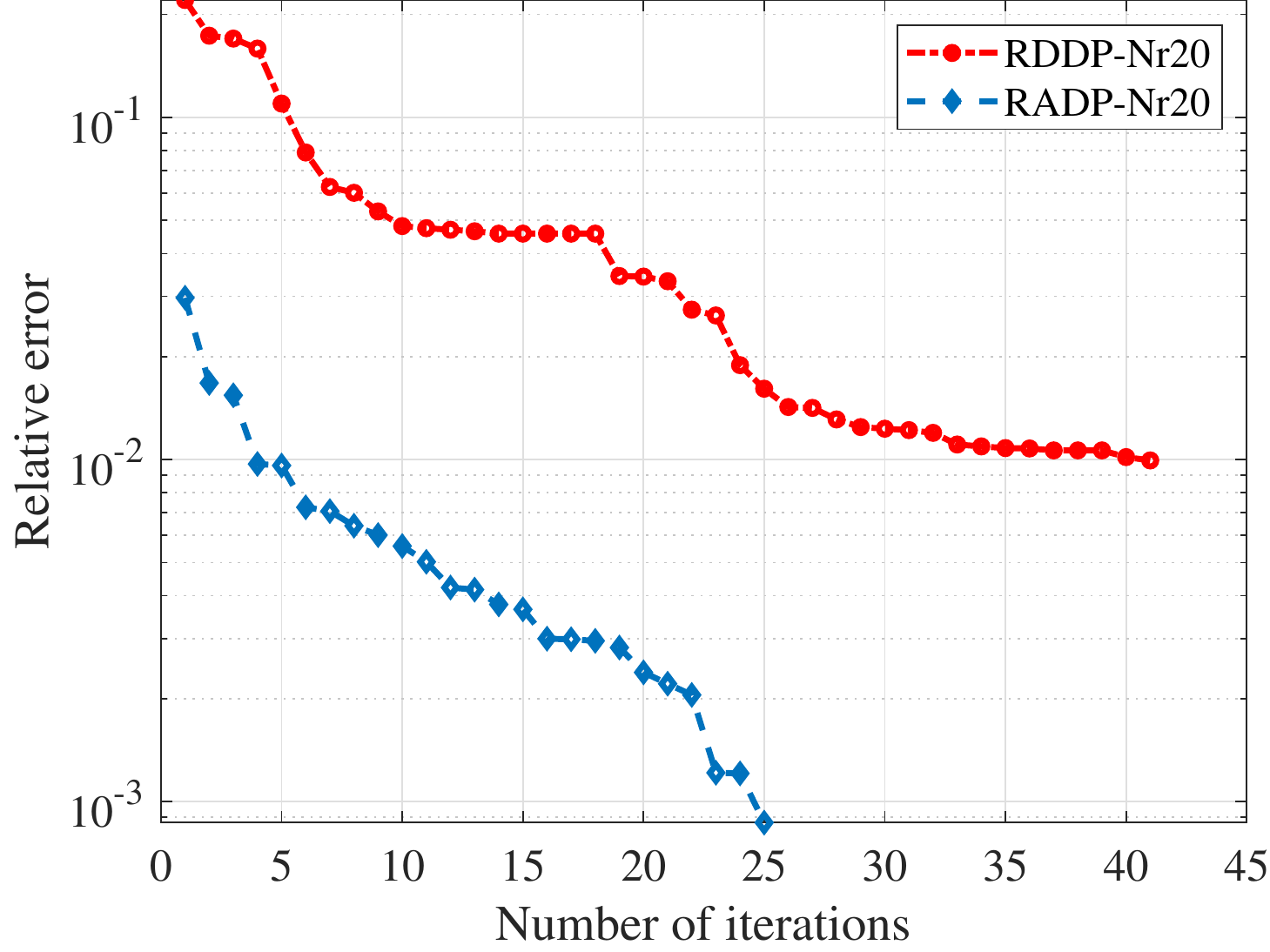}} 
    \subfloat
    {\label{Nr90_fig}\includegraphics[width=0.32\textwidth]{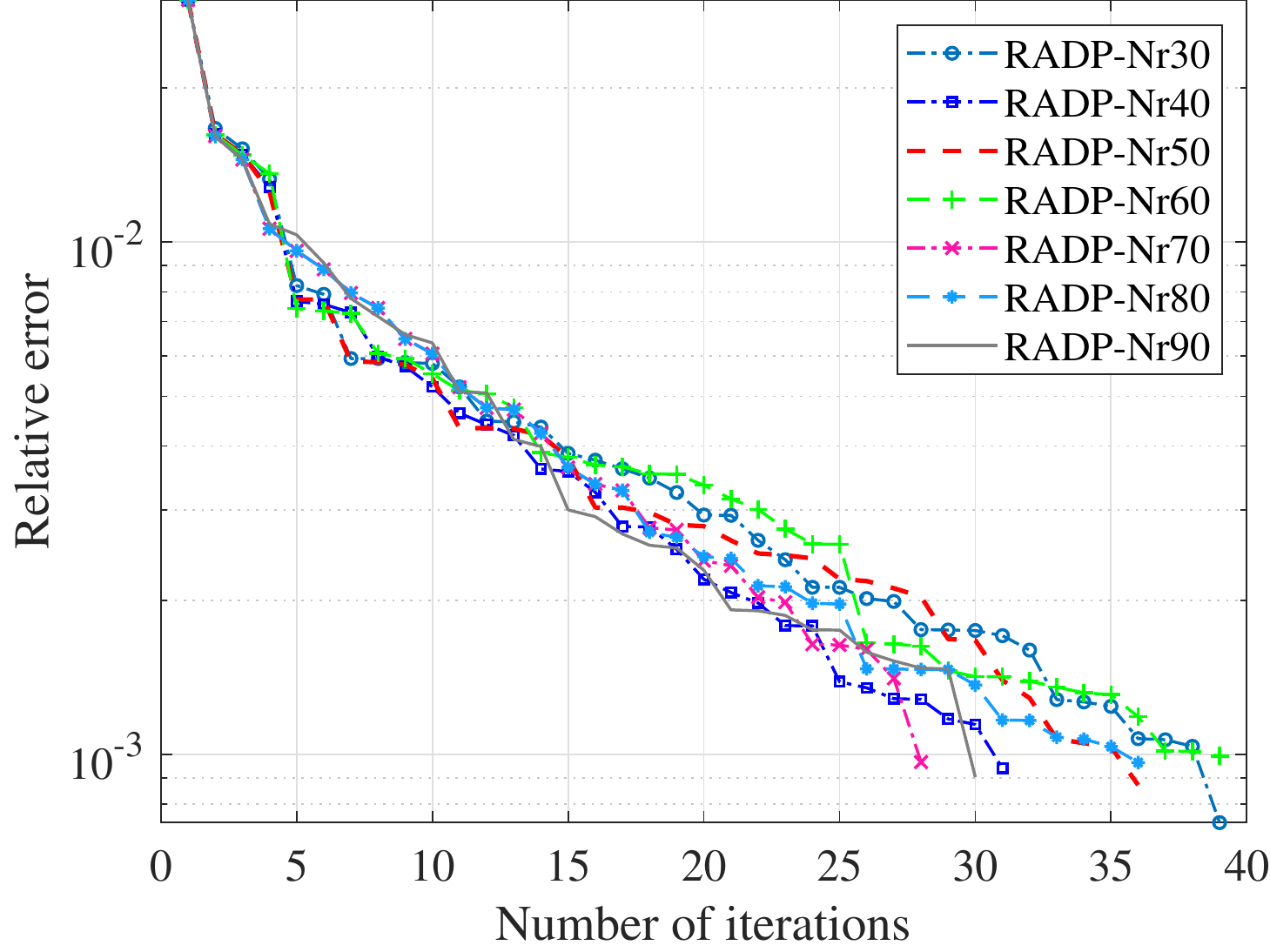}}
    \caption{{Iterative updates of the relative errors between the upper and lower bounds for 118-bus.}}
    \label{fig:RADP/RDDP}
\end{figure*}

We test all the methods for solving the multistage robust UC in 24 stages on the test cases. Apart from the proposed RADP, the comparison methods include the multistage affine decision rule (MAR) method and the robust feasible region (RFR) method. The MAR method models the dispatch decisions as the simplified rule $\boldsymbol{w}_t+\boldsymbol{W}_t\sum \xii_t$, where $\boldsymbol{w}_t$ and $\boldsymbol{W}_t$ are the multistage affine parameters needed to be determined. The difference between the RFR and RADP is that RFR solves the $F^{\rm{uc}}$ problem iteratively until finding the suboptimal solution without approaching the cost-to-go functions in each ED stage. 
Here the comparisons between RADP and RFR are to show the economical performance with/without the cost-to-go functions under the randomized Monte-Carlo realizations. 
The stop criteria for RADP is based on the consecutive change of the objective value $F^{\rm{uc}}$, while the converged criterion for the cost-to-go approximations is based on the relative errors between the upper and lower bounds. The iterative MAR terminates when the obtained commitment solution and the decision rule are feasible for all the constraints. The largest number of iterations in all the tests is set to 100.

To demonstrate the performance of RADP in approximating the cost-to-go functions when given the commitment decisions, the relative error $(\bar{\mathcal{F}}_1-\FF_1)/\FF_1$ for the 118-bus test case is plotted in Fig. \ref{fig:RADP/RDDP}. 
The 'Nr10' $\sim$ 'Nr90' in the legend means the dimension of the uncertainty in the test. 
Meanwhile, the relative error updates of RDDP for worst-case cost-to-go approximations have also been plotted for comparison, where the upper-bound problems \eqref{upper primal problem} are solved by MILPs \cite{georghiou2019robust}.
Note that all the tests of RDDP/RADP start with initializing the candidate points used in upper-bound problems based on Proposition \ref{upperbound-feasible} for guaranteeing the feasibility of inner approximations.
RDDP can achieve the optimal solutions as well as the optimal worst-case cost-to-go functions for multistage RO problems with all the continuous decision variables.
From Fig. \ref{fig:RADP/RDDP}, both the RDDP and RADP have the finite termination guarantee. 
RDDP needs more iterations than RADP to terminate as RDDP converges to the optimal worst-case cost-to-go functions while comparing to the suboptimal solutions achieved by RADP. 
RDDP encounters computational burdens when the dimension of the uncertainty becomes large.
To better investigate the computational time improvement, we have listed the iteration number (Iter.) and runtime (RT) of RDDP/RADP in Table \ref{table:runtime of RADP/RDDP}. Compared to the RDDP, RADP scales nicely with the dimension of uncertainty as RADP solves the McCormick relaxation based upper-bound bilinear programming problems for obtaining the worst cases. Since the MILP-based RDDP can not converge within 24 hours when $\mathcal{N}_r$ is larger than 30, the runtime of which is not reported.

\begin{table}[!t] 
\caption{Characteristics of the Test Cases}
\label{table:characteristics of cases}
\centering  
\setlength{\tabcolsep}{0.5mm}{
\begin{tabular}{c|c|c|c|c|c}
\hline
Test Case   & Bus No. & Unit No. & Storage No. & Branch No. & Load No.     \\ \hline  
118-bus   & 118 & 54 &10 & 179 &91  \\
2383-bus  &2383 &323 &40  &2896  &1789 \\        
\hline
\end{tabular}}
\end{table}
\begin{table}[!t] 
\caption{{Runtime of RADP/RDDP for 118-bus}}
\label{table:runtime of RADP/RDDP}
\centering  
\begin{tabular}{l|c|c|c|c|c|c|c|c|c}  
\hline
$\mathcal{N}_r$  & 10 & 20 & 30 & 40 & 50 & 60 & 70 & 80 & 90 \\ \hline  
RDDP Iter.       &65  &41     &N/A  &N/A    &N/A   &N/A  &N/A     &N/A         &N/A     \\ 
RDDP RT(s)     &615.5   &1.2$*10^5$    &N/A &N/A  &N/A      &N/A    &N/A     &N/A           &N/A   \\
RADP Iter.       &30  &25     &39  &31    &36   &39  &28     &36         &30   \\ 
RADP RT(s)     &56.2  &52.5     &134.4  &117.1   &148.2  &183.4  &143.2    &210.5        &186.7   \\ 
\hline
\end{tabular}
\end{table}

\begin{table}[!t] 
\caption{{Runtime of RADP/MAR}}
\label{table:time}
\centering  
\begin{tabular}{p{1.32cm}<{\centering}|p{1.32cm}<{\centering}|p{1.32cm}<{\centering}|p{1.32cm}<{\centering}|p{1.32cm}<{\centering}} 
\hline
Test Case    & MAR & MAR & RADP & RADP    \\
Bus-$\mathcal{N}_r$    &Iter.  & Time(s) &Iter.  & Time(s)  \\ \hline  
118-10     &8      &2.49*$10^3$    &2      &12.75      \\            
118-20     &29     &6.63*$10^3$    &3      &51.35       \\ 
118-30     &9      &1.11*$10^3$    &3      &54.78       \\
118-40     &20     &2.34*$10^3$    &3      &102.39      \\
118-50     &27     &3.43*$10^3$    &5      &199.88      \\
118-60     &17     &2.72*$10^3$    &7      &605.24      \\ 
118-70     &31     &4.02*$10^3$    &6      &404.58      \\
118-80     &27     &3.04*$10^3$    &8      &653.28      \\
118-90     &16     &1.44*$10^3$    &3      &207.28      \\
2383-30    &100    &2.69*$10^5$    &4      &536.23      \\
2383-60    &N/A    &N/A            &4      &532.46          \\
2383-90    &N/A    &N/A            &8      &1.64*$10^3$     \\
2383-120   &N/A    &N/A            &7      &1.78*$10^3$     \\
2383-150   &N/A    &N/A            &10     &3.34*$10^3$     \\
2383-180   &N/A    &N/A            &9      &7.21*$10^3$     \\
2383-210   &N/A    &N/A            &10     &6.12*$10^3$     \\
2383-240   &N/A    &N/A            &5      &5.89*$10^3$     \\
2383-270   &N/A    &N/A            &8      &1.23*$10^4$     \\
\hline
\end{tabular}
\end{table}
\begin{table}[!t] 
\caption{{Operating Cost of MAR/RFR/RADP}}
\label{table:cost}
\centering  
\begin{tabular}{p{1.6cm}<{\centering}|p{1.6cm}<{\centering}|p{1.6cm}<{\centering}|p{1.6cm}<{\centering}} 
\hline
Test Case    & MAR & RFR & RADP    \\
Bus-$\mathcal{N}_r$    & Cost(\$) & Cost(\$) & Cost(\$)   \\ \hline  
118-10     &1.64*$10^6$   &1.63*$10^6$       &1.62*$10^6$    \\            
118-20     &1.64*$10^6$   &1.63*$10^6$       &1.61*$10^6$    \\ 
118-30     &1.64*$10^6$   &1.63*$10^6$       &1.61*$10^6$     \\
118-40     &1.64*$10^6$   &1.63*$10^6$       &1.61*$10^6$     \\
118-50     &1.65*$10^6$   &1.62*$10^6$       &1.60*$10^6$     \\
118-60     &1.65*$10^6$   &1.62*$10^6$       &1.60*$10^6$     \\ 
118-70     &1.65*$10^6$   &1.63*$10^6$       &1.60*$10^6$      \\
118-80     &1.63*$10^6$   &1.63*$10^6$       &1.60*$10^6$      \\
118-90     &1.66*$10^6$   &1.62*$10^6$       &1.60*$10^6$      \\
2383-30    &4.49*$10^7$   &4.48*$10^7$       &4.45*$10^7$      \\
2383-60    &N/A           &4.48*$10^7$       &4.46*$10^7$       \\
2383-90    &N/A           &4.48*$10^7$       &4.45*$10^7$       \\
2383-120   &N/A           &4.47*$10^7$       &4.45*$10^7$       \\
2383-150   &N/A           &4.48*$10^7$       &4.44*$10^7$       \\
2383-180   &N/A           &4.48*$10^7$       &4.45*$10^7$       \\
2383-210   &N/A           &4.48*$10^7$       &4.45*$10^7$       \\
2383-240   &N/A           &4.47*$10^7$       &4.44*$10^7$       \\
2383-270   &N/A           &4.48*$10^7$       &4.45*$10^7$       \\
\hline
\end{tabular}
\end{table}

To better investigate the computational improvement for solving multistage UC problems, we further compare the proposed RADP with MAR on the 118-bus and 2383-bus test cases, under different settings for the dimension of the uncertainty. 
The running time and the iteration number of termination are listed for the test cases in Table \ref{table:time}. 
RADP achieves computational gains when compared to the running time and the iterations of MAR. 
RADP can still obtain the suboptimal solutions for large-scale systems as in 2383-bus test cases, while MAR cannot find solutions in 24 hours.
The main reasons for the computational improvement of RADP are two-fold: 1) the nonanticipativity constrained $F^{\rm{uc}}$ in \eqref{nonancipative UC-stage} in each iteration can be solved faster than the master problem of MAR for obtaining commitment solutions and affine rules, 2) RADP terminates at fewer iteration numbers than MAR as RADP can generate more scenarios in $\mathcal{U}_t$ at each iteration for obtaining the converged cost-to-go functions in the inner loop.

We further evaluate the economical performance of RADP compared with MAR/RFR on the 118-bus and 2383-bus test cases with different $\mathcal{N}_r$ settings. The operating costs are listed in Table \ref{table:cost}, where the operating cost  in each test is computed by averaging over 200 Monte-Carlo realizations. 
In each Monte-Carlo run, MAR calculates the ED solutions using the obtained affine rules based on the specific realizations. 
RADP obtains the solutions at every ED stage by solving $\underaccent{\bar}{F}_{t}(\y_{t-1}, \xii_{t})$ in \eqref{stage t lower-bound problem} with the converged cost-to-go functions. 
The operating cost of RFR is calculated stage by stage based on the obtained robust feasible region limits of the dispatch variables without using the converged cost-to-go functions.
The average operating costs of MAR are larger than RADP and RFR in most of the test cases, which implies that the preselected affine rules may result in uneconomic ED solutions compared to the robust feasible region based methods. Furthermore, the average operating costs of RADP are less than RFR, which means the converged cost-to-go functions can help to achieve economical improvement.  

\section{Conclusions}
\label{section5}
This paper presented a robust dynamic programming based framework for solving the multistage robust UC problem, in order to promote the economical and computational performance for large-scale power system decision-making. 
The RDDP scheme can achieve finite convergence of the optimal worst-case cost-to-go functions with all the continuous variables, while it suffers from high computational complexity in upper-bounding the value functions. 
To tackle the issue, we propose the RADP scheme to solve the multistage robust UC problem with both discrete and continuous decisions. 
We initialize the upper bounds with limited candidate points for guaranteeing the feasibility of inner approximations. 
Furthermore, the worst cases are updated based on the McCormick relaxed upper-bound bilinear programming to accelerate the computational performance. The worst cases are used for obtaining both the cost-to-go functions and the nonanticipativity constrained commitment solutions, combining the dual and primal variable updates for the multistage robust UC problems. 
For the RADP method, the finite termination is guaranteed based on the analysis of the finite upper/lower bound generations and the finite scenarios added in constraint generation for obtaining the commitment decisions. 
Extensive numerical comparisons have demonstrated the improved economical performance of the proposed RADP over MAR/RFR, while greatly reducing the computation time. 

Interesting future research directions open up, including the extensions for the quadratic or nonconvex cost-to-go approximation based multistage robust optimization problems. Moreover, we are interested to pursue decentralized dual dynamic programming for large-scale decision-making problems under an uncertain environment.







\section*{Acknowledgment}
This work is supported by National Key R\&D Program of China (2022YFA1004601), National Natural Science Foundation of China (62103323, 11991023, 11991020, 61902308, 62192750), Initiative Postdocs Supporting Program (BX20200270, BX20190275), China Postdoctoral Science Foundation (2021M692565, 2019M663723), and the Fundamental Research Funds for the Central Universities under grant (xhj032021013, xxj022019016).

\ifCLASSOPTIONcaptionsoff
\fi
\bibliographystyle{IEEEtran}   
\bibliography{ref}

\begin{thebibliography}{10}
\providecommand{\url}[1]{#1}
\csname url@samestyle\endcsname
\providecommand{\newblock}{\relax}
\providecommand{\bibinfo}[2]{#2}
\providecommand{\BIBentrySTDinterwordspacing}{\spaceskip=0pt\relax}
\providecommand{\BIBentryALTinterwordstretchfactor}{4}
\providecommand{\BIBentryALTinterwordspacing}{\spaceskip=\fontdimen2\font plus
\BIBentryALTinterwordstretchfactor\fontdimen3\font minus
  \fontdimen4\font\relax}
\providecommand{\BIBforeignlanguage}[2]{{%
\expandafter\ifx\csname l@#1\endcsname\relax
\typeout{** WARNING: IEEEtran.bst: No hyphenation pattern has been}%
\typeout{** loaded for the language `#1'. Using the pattern for}%
\typeout{** the default language instead.}%
\else
\language=\csname l@#1\endcsname
\fi
#2}}
\providecommand{\BIBdecl}{\relax}
\BIBdecl

\bibitem{lorca2016multistageTPS}
A.~Lorca and X.~A. Sun, ``Multistage robust unit commitment with dynamic
  uncertainty sets and energy storage,'' \emph{IEEE Transactions on Power
  Systems}, vol.~32, no.~3, pp. 1678--1688, 2016.

\bibitem{bertsimas2012adaptive}
D.~Bertsimas, E.~Litvinov, X.~A. Sun, J.~Zhao, and T.~Zheng, ``Adaptive robust
  optimization for the security constrained unit commitment problem,''
  \emph{IEEE Transactions on Power Systems}, vol.~28, no.~1, pp. 52--63, 2012.

\bibitem{lorca2016multistageRO}
A.~Lorca, X.~A. Sun, E.~Litvinov, and T.~Zheng, ``Multistage adaptive robust
  optimization for the unit commitment problem,'' \emph{Operations Research},
  vol.~64, no.~1, pp. 32--51, 2016.

\bibitem{wang2008security}
J.~Wang, M.~Shahidehpour, and Z.~Li, ``Security-constrained unit commitment
  with volatile wind power generation,'' \emph{IEEE Transactions on Power
  Systems}, vol.~23, no.~3, pp. 1319--1327, 2008.

\bibitem{wu2007stochastic}
L.~Wu, M.~Shahidehpour, and T.~Li, ``Stochastic security-constrained unit
  commitment,'' \emph{IEEE Transactions on power systems}, vol.~22, no.~2, pp.
  800--811, 2007.

\bibitem{zheng2014stochastic}
Q.~P. Zheng, J.~Wang, and A.~L. Liu, ``Stochastic optimization for unit
  commitment—a review,'' \emph{IEEE Transactions on Power Systems}, vol.~30,
  no.~4, pp. 1913--1924, 2014.

\bibitem{wang2011chance}
Q.~Wang, Y.~Guan, and J.~Wang, ``A chance-constrained two-stage stochastic
  program for unit commitment with uncertain wind power output,'' \emph{IEEE
  transactions on power systems}, vol.~27, no.~1, pp. 206--215, 2011.

\bibitem{wu2014chance}
H.~Wu, M.~Shahidehpour, Z.~Li, and W.~Tian, ``Chance-constrained day-ahead
  scheduling in stochastic power system operation,'' \emph{IEEE Transactions on
  Power Systems}, vol.~29, no.~4, pp. 1583--1591, 2014.

\bibitem{jiang2011robust}
R.~Jiang, J.~Wang, and Y.~Guan, ``Robust unit commitment with wind power and
  pumped storage hydro,'' \emph{IEEE Transactions on Power Systems}, vol.~27,
  no.~2, pp. 800--810, 2011.

\bibitem{zeng2013solving}
B.~Zeng and L.~Zhao, ``Solving two-stage robust optimization problems using a
  column-and-constraint generation method,'' \emph{Operations Research
  Letters}, vol.~41, no.~5, pp. 457--461, 2013.

\bibitem{jabr2014robust}
R.~A. Jabr, S.~Karaki, and J.~A. Korbane, ``Robust multi-period opf with
  storage and renewables,'' \emph{IEEE Transactions on Power Systems}, vol.~30,
  no.~5, pp. 2790--2799, 2014.

\bibitem{li2015adjustable}
Z.~Li, W.~Wu, B.~Zhang, and B.~Wang, ``Adjustable robust real-time power
  dispatch with large-scale wind power integration,'' \emph{IEEE transactions
  on sustainable energy}, vol.~6, no.~2, pp. 357--368, 2015.

\bibitem{attarha2019affinely}
A.~Attarha, P.~Scott, and S.~Thi{\'e}baux, ``Affinely adjustable robust admm
  for residential der coordination in distribution networks,'' \emph{IEEE
  Transactions on Smart Grid}, vol.~11, no.~2, pp. 1620--1629, 2019.

\bibitem{he2016robust}
C.~He, L.~Wu, T.~Liu, and M.~Shahidehpour, ``Robust co-optimization scheduling
  of electricity and natural gas systems via admm,'' \emph{IEEE Transactions on
  Sustainable Energy}, vol.~8, no.~2, pp. 658--670, 2016.

\bibitem{zhai2016transmission}
Q.~Zhai, X.~Li, X.~Lei, and X.~Guan, ``Transmission constrained uc with wind
  power: An all-scenario-feasible milp formulation with strong
  nonanticipativity,'' \emph{IEEE Transactions on Power Systems}, vol.~32,
  no.~3, pp. 1805--1817, 2016.

\bibitem{cobos2018robust}
N.~G. Cobos, J.~M. Arroyo, N.~Alguacil, and J.~Wang, ``Robust energy and
  reserve scheduling considering bulk energy storage units and wind
  uncertainty,'' \emph{IEEE Transactions on Power Systems}, vol.~33, no.~5, pp.
  5206--5216, 2018.

\bibitem{cobos2018robustTSE}
N.~G. Cobos, J.~M. Arroyo, N.~Alguacil, and A.~Street, ``Robust energy and
  reserve scheduling under wind uncertainty considering fast-acting
  generators,'' \emph{IEEE Transactions on Sustainable Energy}, vol.~10, no.~4,
  pp. 2142--2151, 2018.

\bibitem{li2019multi}
X.~Li and Q.~Zhai, ``Multi-stage robust transmission constrained unit
  commitment: A decomposition framework with implicit decision rules,''
  \emph{International Journal of Electrical Power \& Energy Systems}, vol. 108,
  pp. 372--381, 2019.

\bibitem{zhou2020multistage}
Y.~Zhou, Q.~Zhai, and L.~Wu, ``Multistage transmission-constrained unit
  commitment with renewable energy and energy storage: implicit and explicit
  decision methods,'' \emph{IEEE Transactions on Sustainable Energy}, vol.~12,
  no.~2, pp. 1032--1043, 2020.

\bibitem{georghiou2019robust}
A.~Georghiou, A.~Tsoukalas, and W.~Wiesemann, ``Robust dual dynamic
  programming,'' \emph{Operations Research}, vol.~67, no.~3, pp. 813--830,
  2019.

\bibitem{shapiro2011analysis}
A.~Shapiro, ``Analysis of stochastic dual dynamic programming method,''
  \emph{European Journal of Operational Research}, vol. 209, no.~1, pp. 63--72,
  2011.

\bibitem{philpott2013solving}
A.~Philpott, V.~de~Matos, and E.~Finardi, ``On solving multistage stochastic
  programs with coherent risk measures,'' \emph{Operations Research}, vol.~61,
  no.~4, pp. 957--970, 2013.

\bibitem{baucke2017deterministic}
R.~Baucke, A.~Downward, and G.~Zakeri, ``A deterministic algorithm for solving
  multistage stochastic programming problems,'' \emph{Optimization Online}, pp.
  1--25, 2017.

\bibitem{moazeni2018risk}
S.~Moazeni, A.~H. Miragha, and B.~Defourny, ``A risk-averse stochastic dynamic
  programming approach to energy hub optimal dispatch,'' \emph{IEEE
  Transactions on Power Systems}, vol.~34, no.~3, pp. 2169--2178, 2018.

\bibitem{shuai2018stochastic}
H.~Shuai, J.~Fang, X.~Ai, Y.~Tang, J.~Wen, and H.~He, ``Stochastic optimization
  of economic dispatch for microgrid based on approximate dynamic
  programming,'' \emph{IEEE Transactions on Smart Grid}, vol.~10, no.~3, pp.
  2440--2452, 2018.

\bibitem{xu2013kernel}
X.~Xu, C.~Lian, L.~Zuo, and H.~He, ``Kernel-based approximate dynamic
  programming for real-time online learning control: An experimental study,''
  \emph{IEEE Transactions on Control Systems Technology}, vol.~22, no.~1, pp.
  146--156, 2013.

\bibitem{zeng2018dynamic}
P.~Zeng, H.~Li, H.~He, and S.~Li, ``Dynamic energy management of a microgrid
  using approximate dynamic programming and deep recurrent neural network
  learning,'' \emph{IEEE Transactions on Smart Grid}, vol.~10, no.~4, pp.
  4435--4445, 2018.

\bibitem{lu2019multi}
R.~Lu, T.~Ding, B.~Qin, J.~Ma, X.~Fang, and Z.~Dong, ``Multi-stage stochastic
  programming to joint economic dispatch for energy and reserve with uncertain
  renewable energy,'' \emph{IEEE Transactions on Sustainable Energy}, vol.~11,
  no.~3, pp. 1140--1151, 2019.

\bibitem{papavasiliou2017application}
A.~Papavasiliou, Y.~Mou, L.~Cambier, and D.~Scieur, ``Application of stochastic
  dual dynamic programming to the real-time dispatch of storage under renewable
  supply uncertainty,'' \emph{IEEE Transactions on Sustainable Energy}, vol.~9,
  no.~2, pp. 547--558, 2017.

\bibitem{shi2020enhancing}
Y.~Shi, S.~Dong, C.~Guo, Z.~Chen, and L.~Wang, ``Enhancing the flexibility of
  storage integrated power system by multi-stage robust dispatch,'' \emph{IEEE
  Transactions on Power Systems}, vol.~36, no.~3, pp. 2314--2322, 2020.

\bibitem{lan2022fast}
Y.~Lan, Q.~Zhai, X.~Liu, and X.~Guan, ``Fast stochastic dual dynamic
  programming for economic dispatch in distribution systems,'' \emph{IEEE
  Transactions on Power Systems}, 2022.

\bibitem{bhattacharya2016managing}
A.~Bhattacharya, J.~P. Kharoufeh, and B.~Zeng, ``Managing energy storage in
  microgrids: A multistage stochastic programming approach,'' \emph{IEEE
  Transactions on Smart Grid}, vol.~9, no.~1, pp. 483--496, 2016.

\bibitem{shi2019multistage}
Z.~Shi, H.~Liang, S.~Huang, and V.~Dinavahi, ``Multistage robust energy
  management for microgrids considering uncertainty,'' \emph{IET Generation,
  Transmission \& Distribution}, vol.~13, no.~10, pp. 1906--1913, 2019.

\bibitem{zou2019stochastic}
J.~Zou, S.~Ahmed, and X.~A. Sun, ``Stochastic dual dynamic integer
  programming,'' \emph{Mathematical Programming}, vol. 175, no.~1, pp.
  461--502, 2019.

\bibitem{zou2018multistage}
------, ``Multistage stochastic unit commitment using stochastic dual dynamic
  integer programming,'' \emph{IEEE transactions on Power Systems}, vol.~34,
  no.~3, pp. 1814--1823, 2018.

\bibitem{hjelmeland2018nonconvex}
M.~N. Hjelmeland, J.~Zou, A.~Helseth, and S.~Ahmed, ``Nonconvex medium-term
  hydropower scheduling by stochastic dual dynamic integer programming,''
  \emph{IEEE Transactions on Sustainable Energy}, vol.~10, no.~1, pp. 481--490,
  2018.

\bibitem{ding2020multi}
T.~Ding, M.~Qu, C.~Huang, Z.~Wang, P.~Du, and M.~Shahidehpour, ``Multi-period
  active distribution network planning using multi-stage stochastic programming
  and nested decomposition by sddip,'' \emph{IEEE Transactions on Power
  Systems}, vol.~36, no.~3, pp. 2281--2292, 2020.

\bibitem{xiong2022multi}
H.~Xiong, Y.~Shi, Z.~Chen, C.~Guo, and Y.~Ding, ``Multi-stage robust dynamic
  unit commitment based on pre-extended-fast robust dual dynamic programming,''
  \emph{IEEE Transactions on Power Systems}, 2022.

\bibitem{lee2013modeling}
C.~Lee, C.~Liu, S.~Mehrotra, and M.~Shahidehpour, ``Modeling transmission line
  constraints in two-stage robust unit commitment problem,'' \emph{IEEE
  Transactions on Power Systems}, vol.~29, no.~3, pp. 1221--1231, 2013.

\bibitem{mccormick1976computability}
G.~P. McCormick, ``Computability of global solutions to factorable nonconvex
  programs: Part i—convex underestimating problems,'' \emph{Mathematical
  programming}, vol.~10, no.~1, pp. 147--175, 1976.

\bibitem{ben2017convex}
W.~Ben-Ameur, A.~Ouorou, and G.~Wang, ``Convex and concave envelopes: revisited
  and new perspectives,'' \emph{Operations Research Letters}, vol.~45, no.~5,
  pp. 421--426, 2017.

\bibitem{deng2021optimal}
L.~Deng, H.~Sun, B.~Li, Y.~Sun, T.~Yang, and X.~Zhang, ``Optimal operation of
  integrated heat and electricity systems: A tightening mccormick approach,''
  \emph{Engineering}, vol.~7, no.~8, pp. 1076--1086, 2021.

\bibitem{zimmerman2010matpower}
R.~D. Zimmerman, C.~E. Murillo-S{\'a}nchez, and R.~J. Thomas, ``Matpower:
  Steady-state operations, planning, and analysis tools for power systems
  research and education,'' \emph{IEEE Transactions on power systems}, vol.~26,
  no.~1, pp. 12--19, 2010.

\end{thebibliography}

\end{document}